\documentclass[a4paper]{amsart}
\usepackage{latexsym}\usepackage{ifthen}\usepackage[leqno]{amsmath}
\usepackage{enumerate}\usepackage{calc}\usepackage{hyphenat}
\usepackage{amstext,amsbsy,amsopn,amsthm,amsgen}
\usepackage{amsfonts,amscd,amsxtra,upref}
\usepackage{graphics}
\usepackage{epstopdf}
\usepackage{dsfont}
\usepackage{hyperref}
\usepackage{mathrsfs}\usepackage{euscript}\usepackage{amssymb}
\swapnumbers
\theoremstyle{plain}

\makeatother
\newtheorem{Thm}{Theorem}[section]
\newtheorem{Lem}[Thm]{Lemma}
\newtheorem{Cor}[Thm]{Corollary}
\newtheorem{Pro}[Thm]{Proposition}

\theoremstyle{definition}
\newtheorem{Def}[Thm]{Definition}

\newtheorem{Prb}[Thm]{Problem}
\theoremstyle{remark}
\newtheorem{Rem}[Thm]{Remark}
\numberwithin{equation}{section}

\newcommand{\ITE}[3]{\ifthenelse{#1}{#2}{#3}}\newcommand{\ITEE}[4][]{\ITE{\equal{#2}{#3}}{#4}{#1}}
 \ITE{\isundefined{\texorpdfstring}}{\newcommand{\texorpdfstring}[2]{#1}}{}

\newcommand{\myData}[1][]{
 \author[D.\ Burek]{Dominik Burek}
 \address{\ITEE{#1}{*}{D.\ Burek{}\\{}}
  Instytut Matematyki\\{}Wydzia\l{} Matematyki i~Informatyki\\{}Uniwersytet Jagiello\'{n}ski\\{}
  ul.\ \L{}ojasiewicza 6\\{}30-348 Krak\'{o}w\\{}Poland}
 \email{dominik.burek@student.uj.edu.pl}
 }

\newcommand{\orb}{\operatorname{orb}}
\newcommand{\age}{\operatorname{age}}

\begin{document}
\title{Hodge numbers of generalised Borcea-Voisin threefolds}
\myData

\begin{abstract}
We shall reproof formulas for the Hodge numbers of Calabi-Yau threefolds of Borcea-Voisin type constructed by A.\ Cattaneo and A.\ Garbagnati, using the orbifold cohomology formula and the orbifold Euler characteristic.
\end{abstract}

\subjclass[]{Primary 14J32; Secondary 14J28, 14J17.}
\keywords{Borcea-Voisin, Calabi-Yau 3-folds, orbifold cohomology, orbifold Euler characteristic.}
\thanks{This work was partially supported by the grant 346300 for IMPAN from the Simons Foundation and the matching 2015-2019 Polish MNiSW fund.}
\maketitle

\section{Introduction} 

One of the many reasons behind the interest in non-symplectic automorphisms of $K3$ surfaces is the mirror symmetry construction of C.\ Borcea (\cite{B}) and C.\ Voisin (\cite{V}). They independently constructed a family of Calabi-Yau threefolds using a non-symplectic involutions of $K3$ surfaces and elliptic curves. Moreover C.\ Voisin gave a construction of explicit mirror maps. 

\begin{Thm}[\cite{B} and \cite{V}]{} Let $E$ be an elliptic curve with an involution $\alpha_{E}$ which does not preserve $\omega_{E}$. Let $S$ be a $K3$ surface with a non-symplectic involution $\alpha_{S}$. Then any crepant resolution of the variety $(E\times  S)/(\alpha_{E}\times \alpha_{S})$ is a Calabi-Yau manifold with 
$$h_{}^{1,1}=11+5N'-N\quad \textrm{and}\quad h_{}^{2,1}=11+5N-N',$$
where $N$ is a number of curves in $S^{\alpha_{S}}$ and $N'$ is a sum of their genera.
\end{Thm}

In \cite{C} A.\ Cattaneo and A.\ Garbagnati generalised Borcea-Voisin construction using purely non-symplectic automorphisms of order 3, 4 and 6. They obtained the following theorem

\begin{Thm}[\cite{C}]{}
Let $S$ be a $K3$ surface admitting a purely non-symplectic automorphism $\alpha_{S}$ of order $n.$ Let $E$ be an elliptic curve
admitting an automorphism $\alpha_{E}$ such that $\alpha_{E}(\omega_{E})=\zeta_{n}\omega_{E}.$ Then $n \in\{2,3,4,6\}$ and
$({S\times E})/(\alpha_{S}\times \alpha_{E}^{n-1})$ is a singular variety which admits a crepant resolution which
is a Calabi-Yau manifold.
\end{Thm}

Any crepant resolution of $({S\times E})/(\alpha_{S}\times \alpha_{E}^{n-1})$ will be called a
\textsl{Calabi-Yau 3-fold of Borcea-Voisin type.}
For all possible orders they computed the Hodge numbers of this varieties and constructed an elliptic fibrations on them. Their computations are more technical and rely on a detailed study of a crepant resolutions of threefolds. \par
In this paper we give shorter computations of the Hodge numbers using orbifold cohomology introducted by W.\ Chen and Y.\ Ruan in \cite{CR} and orbifold Euler characteristic (cf. \cite{R3}). The main advantage of our approach is that the computations are carried out on $S\times E.$

\section{Preliminaries}

A Calabi-Yau manifold $X$ is a complex, smooth, projective $d$-fold $X$ satisfying 
\begin{enumerate}
\item $K_{X}=\mathcal{O}_{X},$
\item $H^{i}\mathcal{O}_{X}=0$ for $0<i<d.$
\end{enumerate}
A non-trivial generator $\omega_{X}$ of $H^{d,0}(X)\simeq \mathbb{C}$ is called a \textit{period} of $X.$ For any automorphism $\alpha_{X}\in \textup{Aut}(X),$ the induced mapping $\alpha_{X}^{*}$ acts on $H^{d,0}(X)$ and $\alpha_{X}^{*}(\omega_{X})=\lambda_{\alpha}\omega_{X},$ for some $\lambda_{\alpha}\in \mathbb{C}^{*}.$\par
In the case of a $K3$ surface $S$ automorphism $\alpha_{S}\in \textup{Aut}(X)$ which preserves a period is called \textit{symplectic.} If $\alpha_{S}$ does not preserve a period then it is called \textit{non-symplectic.} If additionally $\alpha_{S}$ is of finite order $n$ and $\alpha_{S}^{*}(\omega_{S})=\zeta_{n}\omega_{S},$ where $\zeta_{n}$ is a primitive $n$-th root of unity, then it is called \textit{purely non-symplectic.} \par
\par
We have the following characterisation of orders of purely non-symplectic automorphisms of elliptic curves and $K3$ surfaces.

\begin{Thm}[\cite{S}]{}\label{S} Let $E$ be an elliptic curve. If $\alpha_{E}$ is an automorphism of $E$ which does not preserve
the period, then $\alpha_{E}(\omega_{E})=\zeta_{n}\omega_{E}$ where $n = 2,3,4,6.$ \end{Thm}

\begin{Thm}[\cite{INK}]{}\label{INK} Let $S$ be a $K3$ surface and $\alpha_{S}$ be a purely non-symplectic
automorphism of order $n.$ Then $n\leq 66$ and if $n = p$ is a prime number, then
$p\leq 19.$ \end{Thm}

Let $\sigma$ be a purely non-symplectic automorphism of order $p$ of a $K3$ surface $S.$ We will denote by $S^{\sigma}$ the set of fixed points of $\sigma.$ The action of $\sigma$ may be locally linearized and diagonalized at $p\in S^{\sigma}$ (cf. Cartan \cite{Car}), so the possible local actions are 
$$
{
\begin{pmatrix}
   \zeta_{p}^{t+1} & 0 \\
  0 & \zeta_{p}^{p-t}  \\
\end{pmatrix},
}\quad
\textrm{for } t=0,1,\ldots, p-2.
$$
Clearly, if $t=0,$ then $p$ belongs to a smooth curve fixed by $\sigma,$ otherwise $p$ is an isolated point.
We have the following description of the fixed locus of $\sigma.$

\begin{Thm}[\cite{MAS}, Lemma 2.2, p. 5]{}
Let $S$ be a $K3$ surface and let $\alpha_{S}$ be a non-symplectic automorphism of S of order $n$. Then there are three possibilities
\begin{itemize} 
\item  $S^{\alpha_{S}}= \emptyset;$ in this case $n=2,$ 
\item  $S^{\alpha_{S}} = E_{1}\cup E_{2},$ where $E_{1},$ $E_{2}$ are disjoint smooth elliptic curves; in this case $n=2,$
\item  $S^{\alpha_{S}}=C\cup R_{1}\cup R_{2}\cup \ldots R_{k-1}\cup \{p_{1}, p_{2}, \ldots, p_{h}\},$ where $p_{i}$ are isolated fixed points, $R_{i}$ are smooth rational curves and $C$ is the curve with highest genus $g(C).$
\end{itemize}
\end{Thm}
We refer to \cite{MAS} for a proof and more precise description of the fixed locus for particular values of $n$.\par

\section{Orbifold's cohomology}

In \cite{CR} W.\ Chen and Y.\ Ruan introduced a new cohomology theory for orbifolds. We consider varieties $X/G,$ where $X$ is a projective variety and $G$ is a finite group acting on $X$ viewed as orbifold.

\begin{Def} For $G \in GL_{n}(\mathbb{C})$ of order $m,$ let $e^{2\pi i a_{1}}, e^{2\pi i a_{2}}, \ldots, e^{2\pi i a_{n}}$ be eigenvalues of $G$ for some $a_{1}, a_{2}, \ldots, a_{n}\in [0,1)\cap \mathbb{Q}$.  The value of the sum $a_1 + a_2+ \ldots + a_n$ is called the age of $G$ and is denoted by $\age(G).$
\end{Def}

The age of $G$ is an integer if and only if $\det G=1$ i.e. $G \in SL_{n}(\mathbb{C})$.

\begin{Def} For a variety $X/G$ define the Chen-Ruan cohomology by $$H_{\textrm{orb}}^{i,j}(X/G):=\bigoplus_{[g]\in\textrm{Conj}(G)}\left(\bigoplus_{U\in \Lambda(g)} H^{i- \age(g),\; j-\age(g)}(U)\right)^{\textup{C}(g)},$$ where $\textrm{Conj}(G)$ is the set of conjugacy classes of $G$ (we choose a representative $g$ of each conjugacy class), $\textup{C}(g)$ is the centralizer of $g$, $\Lambda(g)$ denotes the set of irreducible connected components of the set fixed by $g\in G$ and $\age(g)$ is the age of the matrix of linearized action of $g$ near a point of $U.$ \par The dimension of $ H_{\textrm{orb}}^{i,j}(X/G)$ will be denoted by $h_{\text{orb}}^{i,j}(X/G).$\end{Def}

\begin{Rem} If the group $G$ is cyclic of a prime order $p$, then we can pick a generator $\alpha$ and the above formula simplifies to 
$$H_{\textrm{orb}}^{i,j}(X/G)=H^{i,j}(X)^{G}\oplus \bigoplus_{U\in \Lambda(\alpha)} \bigoplus_{k =1}^{p-1} H^{i- \age(\alpha^{k}),\; j-\age(\alpha^{k})}(U).$$
\end{Rem}

We have the following theorem.

\begin{Thm}[\cite{TY}, Theorem 1.1, p. 2]{}\label{TY3}
Let $G$ be a finite group acting on an algebraic smooth variety $X$. If there exists a crepant resolution $\widetilde{X/G}$ of variety $X/G,$ then the following equality holds $$h^{i,j}(\widetilde{X/G})=h^{i,j}_{\orb}(X/G).$$
\end{Thm}

\section{Orbifold Euler characteristic}

Let $G$ be a finite group acting on an algebraic variety $X.$ In a similar manner as in the case of Hodge numbers, we can use an orbifold formula to compute the Euler characteristic of a crepant resolution of $X/G$ (for details see \cite{R3}).

\begin{Def} The orbifold Euler characteristic of $X/G$ is defined as $$e_{\orb}(X/G):=\frac{1}{\# G}\sum_{\substack{(g,h)\in G\times G \\ gh=hg}}e(X^{g}\cap X^{h}).$$ \end{Def}

\begin{Thm}[\cite{R3}, Theorem 2, p. 534]{}\label{TR3}
For any finite abelian group $G$ acting on smooth algebraic variety $X$. If there exists a crepant resolution $\widetilde{X/G}$ of variety $X/G,$ then the following equality holds $$e(\widetilde{X/G})=e_{\orb}(X/G).$$ 
\end{Thm}

\section{Computations of Hodge numbers}
\subsection{Order 2}
Let $(S, \alpha_{S})$ be a $K3$ surface admitting a non-symplectic involution $\alpha_{S}$. Consider an elliptic curve $E$ with non-symplectic involution $\alpha_{E}$ (any elliptic curve $E$ admits such an automorphism). Let us denote by $H^{2}(S,\mathbb{C})^{\alpha_{S}}$ the invariant part of cohomology $H^{2}(S,\mathbb{C})$ under $\alpha_{S}$ and by $r$ the dimension $r=\dim H^{2}(S,\mathbb{C})^{\alpha_{S}}.$ We also denote the eigenspace for $-1$ of the induced action $\alpha_{S}^{*}$ on $H^{2}(S,\mathbb{C})$ by $H^{2}(S,\mathbb{C})_{-1}$ and by $m$ the dimension $m=\dim H^{2}(S,\mathbb{C})_{-1}.$ \par
We see that $$H^{2}(S,\mathbb{C})=H^{2}(S,\mathbb{C})^{\alpha_{S}}\oplus H^{2}(S,\mathbb{C})_{-1},$$ hence the Hodge diamonds of the respective eigenspaces have the following forms

\begin{figure}[h]
 \centering 
\begin{minipage}{.4\textwidth}
    \centering
    $H^{i,j}(S,\mathbb{C})^{\alpha_{S}}$ \par\bigskip
    \includegraphics{diamond.1}
\end{minipage}
\begin{minipage}{.4\textwidth}
    \centering
    $H^{i,j}(S,\mathbb{C})_{-1}$ \par\bigskip
    \includegraphics{diamond.2}
\end{minipage}
\end{figure}

The Hodge diamonds of eigenspaces of the induced action of $\alpha_{E}^{*}$ on $H^{1}(E,\mathbb{C})$ have forms

$$
{
\stackrel{\mbox{$H^{i,j}(E,\mathbb{C})^{\alpha_{E}}$}}
{\begin{matrix}
 &   &  &   &  \\
  &   & 1 &   &  \\
  & 0 &   & 0 &  \\
  &   & 1 &   & \\ 
\end{matrix}
}
}
\\ \quad \quad \quad \quad \quad \quad 
{
\stackrel{\mbox{$H^{i,j}(E,\mathbb{C})_{-1}$}}
{\begin{matrix}
 &   &  &   &  \\
  &   & 0 &   &  \\
  & 1 &   & 1 &  \\
  &   & 0 &   & \\ 
\end{matrix}
}
}
$$

By K{\"u}nneth's formula the Hodge diamond of $H^{3}(S\times E, \mathbb{C})^{C_{2}}$ is given by

\begin{figure}[h]
 \centering 
\begin{minipage}{.5\textwidth}
\centering
    \includegraphics{Bigdiam.1}
\end{minipage}
\end{figure}

The local action of $-1$ on curve may be linearized to matrix 
$$
{
\begin{pmatrix}
   1 & 0 & 0\\
  0 & -1 & 0 \\
  0 & 0 & 1 \\
\end{pmatrix}
}
$$ with age equal to 1. Thus by \ref{TR3}
\begin{align*}
h_{}^{i,j}(\widetilde{S\times E/C_{2}})=h_{\orb}^{i,j}(S\times E/C_{2})=h^{i,j}(S\times E)^{C_{2}}\oplus \bigoplus_{U\in \Lambda(-1)} h^{i-1,j-1}(U),
\end{align*}
which gives formulas 
$$h^{1,1}(\widetilde{S\times E/C_{2}})=r + 1 + 4N\;\;\; \textup{and} \;\;\; h_{}^{2,1}(\widetilde{S\times E/C_{2}})=m-1+4N'.$$

Since the quotient $S/\alpha_{S}$ is a smooth surface with Euler characteristic $$e(S/\alpha_{S})=12+N-N',$$ we recover formulas from Thm. 1.1 of \cite{C}.

\subsection{Order 3}
Let $(S, \alpha_{S})$ be a $K3$ surface admitting a purely non-symplectic automorphism $\alpha_{S}$ of order $3.$ Eigenvalues of induced mapping $\alpha_{S}^{*}$ on $H^{2}(S,\mathbb{C})$ belong to $\{1,\zeta_{3}, \zeta_{3}^{2}\}.$ Let us denote by $H^{2}(S,\mathbb{C})_{\zeta_{3}^{i}}$ the eigenspace of the eigenvalue $\zeta_{3}^{i}.$ For $i=1,2$ the dimension of $H^{2}(S,\mathbb{C})_{\zeta_{3}^{i}}$ does not depend on $i$ and will be denoted by $m.$ Moreover let $r$ be a dimension of $H^{2}(S,\mathbb{C})^{\alpha_{S}}$ --- invariant part of $H^{2}(S,\mathbb{C})$ under $\alpha_{S}.$ \par
Consider an elliptic curve $E$ with the Weierstrass equation $y^{2}=x^{3}+1$ together with a non-symplectic automorphism $\alpha_{E}$ of order 3 such that $\alpha_{E}(x,y)=(\zeta_{3}x, y).$ \par
We see that $$H^{2}(S,\mathbb{C})=H^{2}(S,\mathbb{C})^{\alpha_{S}}\oplus H^{2}(S,\mathbb{C})_{\zeta_{3}}\oplus H^{2}(S,\mathbb{C})_{\zeta_{3}^{2}}.$$ Because $\alpha_{S}^{*}\bigr|_{H^{2,0}(S,\mathbb{C})}([\omega])=\zeta_{3}[\omega]$ for any $[\omega]\in H^{2,0}(S),$ we get $H^{2,0}(S)\subset H(S,\mathbb{C})_{\zeta_{3}}.$\newline
The complex conjugation yields $H^{0,2}(S)\subset H^{2}(S,\mathbb{C})_{\zeta_{3}^{2}}.$ Finally $$\alpha_{S}^{*}([\omega_{S}\wedge \omega_{E}])=\bar{\zeta_{3}}\zeta_{3}[\omega_{S}\wedge \omega_{E}]=[\omega_{S}\wedge \omega_{E}],$$ hence the Hodge diamonds of the respective eigenspaces have the following forms

\begin{figure}[h]
 \centering 
\begin{minipage}{.3\textwidth}
    \centering
    $H^{i,j}(S,\mathbb{C})^{\alpha_{S}}$ \par\bigskip
    \includegraphics{Order3.1}
\end{minipage}
\begin{minipage}{.3\textwidth}
    \centering
    $H^{i,j}(S,\mathbb{C})_{\zeta_{3}}$ \par\bigskip
    \includegraphics{Order3.2}
\end{minipage}
\begin{minipage}{.3\textwidth}
    \centering
    $H^{i,j}(S,\mathbb{C})_{\zeta_{3}^{2}}$ \par\bigskip
    \includegraphics{Order3.3}
\end{minipage}
\end{figure}

Similar analysis gives the Hodge diamonds of eigenspaces of action on the Hodge groups. 

$$
{
\stackrel{\mbox{$H^{i,j}(E,\mathbb{C})^{\alpha_{E}^{2}}$}}
{\begin{matrix}
 &   &  &   &  \\
  &   & 1 &   &  \\
  & 0 &   & 0 &  \\
  &   & 1 &   & \\ 
\end{matrix}
}
}
\ \ \ \quad \quad \quad \quad \quad
{
\stackrel{\mbox{$H^{i,j}(E)_{\zeta_{3}}$}}
{\begin{matrix}
 &   &  &   &  \\
  &   & 0 &   &  \\
  & 1 &   & 0 &  \\
  &   & 0 &   & \\ 
\end{matrix}
}
}
\ \ \ \quad \quad \quad \quad \quad
{
\stackrel{\mbox{$H^{i,j}(E)_{\zeta_{3}^{2}}$}}
{\begin{matrix}
 &   &  &   &  \\
  &   & 0 &   &  \\
  & 0 &   & 1 &  \\
  &   & 0 &   & \\ 
\end{matrix}
}
}
$$

By K{\"u}nneth's formula the Hodge diamond of the invariant part of $H^{3}(S\times E, \mathbb{C})$ has the same form as in the case of order 2.


We denote the automorphism $\alpha_{S}\times \alpha_{E}^{2}$ by $\alpha$. Let us now consider possible actions of elements of $\langle\alpha\rangle\simeq C_{3}$ on $S$ and $E.$
\par
\textsl{The action of $\alpha$.} The action of the automorphism $\alpha^{}$ on $E$ is given by $$E \ni (x,y)\mapsto (\zeta_{3}^{2}x, y) \in E,$$ hence it has three fixed points. Locally the action of $\alpha$ on components of the fixed locus can be diagonalised to
$$
{
\begin{pmatrix}
  1 & 0 & 0\\
  0 & \zeta_{3} & 0 \\
  0 & 0 & \zeta_{3}^{2} \\
\end{pmatrix}
}
\ \quad \quad \text{and} \quad \quad
{
\begin{pmatrix}
  \zeta_{3}^{2} & 0 & 0 \\
  0 & \zeta_{3}^{2} & 0 \\
  0 & 0 & \zeta_{3}^{2} \\
\end{pmatrix}.
}
$$

It follows that ages are equal to 1 and 2 respectively.
\par
\textsl{The action of $\alpha^{2}$.} Analogously we gets possible diagonalised matrices
$$
{
\begin{pmatrix}
   \zeta_{3} & 0 & 0\\
  0 & \zeta_{3} & 0 \\
  0 & 0 & \zeta_{3} \\
\end{pmatrix}
}
\ \quad \quad \text{and} \quad \quad
{
\begin{pmatrix}
  1 & 0 & 0 \\
  0 & \zeta_{3}^{2} & 0 \\
  0 & 0 & \zeta_{3} \\
\end{pmatrix},
}
$$
with ages 1 and 1.
\par

Thus decomposing $(S\times E)^{\alpha}=\mathcal{C}\cup \mathcal{R}\cup \mathcal{P}$, where 
\begin{enumerate}[label=--]
\item[] $\mathcal{C}:=\{\textup{3 curves with highest genus $g(C)$} \},$
\item[] $\mathcal{R}:=\{\textup{$3k-3$ rational curves}\},$
\item[] $\mathcal{P}:=\{\textup{$3n$ isolated points}\},$
\end{enumerate}
the orbifold formula implies that
\begin{align*}
&H_{\orb}^{i,j}(S\times E/C_{3})=H^{i,j}(S\times E)^{C_{3}}\oplus \bigoplus_{U\in \Lambda(\zeta_{3})}\bigoplus_{\iota =1}^{2} H^{i- \age(\alpha^{\iota}),\; j-\age(\alpha^{\iota})}(U)= \\
&= H^{i,j}(S\times E)^{C_{3}}\oplus \left(\bigoplus_{U\in \mathcal{C}}H^{i-1,\; j-1}(U)\oplus H^{i-1,\; j-1}(U)\right)\oplus \\ &\oplus\left(\bigoplus_{U\in \mathcal{R}}H^{i-1,\; j-1}(U)\oplus H^{i-1,\; j-1}(U)\right)\oplus \left(\bigoplus_{U\in \mathcal{P}}H^{i-1,\; j-1}(U)\right).
\end{align*}
Therefore by \ref{TY3}
\begin{align*}
&h^{1,1}(\widetilde{X/C_{3}})=r+1 +  6\cdot 1 + 2\cdot (3k-3)\cdot 1 + 3n\cdot 1=r+1+3n+6k.\\
&h^{1,2}(\widetilde{X/C_{3}})=m-1 + 2\cdot 3\cdot g(C)+(3k-3)\cdot 2\cdot 0 +3n\cdot 0 = m-1+6g(C). \end{align*}
\par
Hence we proved the following theorem:

\begin{Thm} If $S^{\alpha_{S}}$ consists of $k$ curves together with a curve with highest genus $g(C)$ and $n$ isolated points, then for any crepant resolution of the variety $(S\times E)/(\alpha_{S}\times \alpha_{E}^{2})$ the following holds
$$ h_{}^{1,1}=r+1+3n+6k\;\;\; \textup{and}\;\;\; h_{}^{2,1}=m-1+6g(C).$$
\end{Thm}

\subsection{Order 4}
Let $(S, \alpha_{S})$ be a $K3$ surface with purely non-symplectic automorphism $\alpha_{S}$ of order $4.$ Consider an elliptic curve $E$ with the Weierstrass equation $y^{2}=x^{3}+x$ together with a non-symplectic automorphism $\alpha_{E}$ of order 4 such that $$\alpha_{E}(x,y)=(-x, iy).$$ Additionally, suppose that $S^{\alpha_{S}^{2}}$ is not a union of two elliptic curves.

We shall keep the notation of \cite{C}.

\begin{enumerate}[label=--]
\item[] $X=S\times E,$
\item[] $P$ -- the infinity point of $E,$
\item[] $r=\dim H^{2}(S,\mathbb{C})^{\alpha_{S}},$
\item[] $m=\dim H^{2}(S,\mathbb{C})_{\zeta_{6}^{i}}$ for $i\in\{1,2,\ldots, 5\},$
\item[] $N$ -- number of curves which are fixed by $\alpha_{S}^{2},$
\item[] $k$ -- number of curves which are fixed by $\alpha_{S}^{}$ (curves of the first type).
\item[]$b$ -- number of curves which are fixed by $\alpha_{S}^{2}$ and are invariant by $\alpha_{S}^{}$ (curves of the second type),
\item[] $a$ -- number of pairs $(A, A')$ of curves which are fixed by $\alpha_{S}^{2}$ and $\alpha_{S}(A)=A'$ (curves of the third type),
\item[] $D$ -- the curve of the highest genus in $S^{\alpha_{S}^{2}},$
\item[] $n_{1}$ -- number of curves which are fixed by $\alpha_{S}$ not laying on the curve $D,$ 
\item[] $n_{2}$ -- number of curves which are fixed by $\alpha_{S}$ laying on the curve $D.$ 
\end{enumerate}


For the same reasons as in the previous cases $$h^{1,1}_{\orb}(X)^{C_{4}}=r+1\;\;\; \textrm{and} \;\;\; h^{2,1}_{\orb}(X)^{C_{4}}=m-1.$$ 

We denote an automorphism $\alpha_{S}\times \alpha_{E}^{3}$ by $\alpha$. Furthermore for any $g\in C_{4}$ let $$M_{g}:=\bigoplus_{U\in \Lambda(g)} H^{1- \age(g),\; 1-\age(g)}(U).$$ We will consider all possible cases.\par
\textsl{The action of $\alpha$ and $\alpha^{3}$.} The action of automorphism $\alpha$ on $E$ is given by $$E \ni (x,y)\mapsto (-x, -iy) \in E,$$ hence it has two fixed points --- $P$ and $(0,0)$. The fixed locus of $\alpha$ on $S$ consists of $k$ curves and $n_{1}+n_{2}$ isolated points. Since locally the action of $\alpha$ on $X$ along the curve can be diagonalised to
$$
{
\begin{pmatrix}
   1 & 0 & 0\\
  0 & i & 0 \\
  0 & 0 & -i \\
\end{pmatrix},
}
$$ we infer that its age equals $1.$ Near a fixed point we have a matrix 

$$
{
\begin{pmatrix}
   -1 & 0 & 0\\
  0 & -i & 0 \\
  0 & 0 & -i \\
\end{pmatrix},
}
$$ hence the age equals 2. 
\par In case of the action of $\alpha^{3}$ we observe that the fixed locus consists of $k$ curves and $n_{1}+n_{2}$ points with ages $1.$ 
We see that the summand of $h^{1,1}_{\textup{orb}}$ from both actions is equal to $2k + 2k+2(n_{1}+n_{2})=4k+2(n_{1}+n_{2})$.  
\newline
\begin{center}
\def\arraystretch{1.4}
\begin{tabular}{c|c|c}
  
  Element of $C_{4}$ & $\alpha$ & $\alpha^{3}$\\
  \hline
  Irreducible components & $2k$ curves, $2n_{1}+2n_{2}$ points & $2k$ curves, $2n_{1}+2n_{2}$\\
  \hline
  The age & curve: 1, point: 2 & curve: 1, point: 1\\
  \hline
  Summand of $h^{1,1}_{\textup{orb}}$ & $2k$ & $2k+2(n_{1}+n_{2})$\\
  
\end{tabular} 
\end{center}
\vspace{2mm}

\textsl{The action of $\alpha^{2}$.} The automorphism $\alpha^{2}$ acts on $E$ as $$E \ni (x,y)\mapsto (x, -y) \in E,$$ hence it has four fixed points --- $P$, $(0,0)$, $(i,0),$ $(-i, 0)$ from which only two are invariant under the action of $\alpha$ and the other two are permuted. After identifying $M_{\zeta_{4}^{2}}$ with the vector space spanned by irreducible components we will find the action of induced map $\alpha^{*}$ on it. \par Because the matrix of the action of ${\alpha_{E}^{2}}^{*}$ on $M_{\zeta_{4}^{2}}$ is
$$
{
\begin{pmatrix}
   1 & 0 & 0 & 0\\
  0 & 1 & 0 & 0\\
  0 & 0 & 0 & 1\\
  0 & 0 & 1 & 0\\
\end{pmatrix}
}
\sim 
{
\begin{pmatrix}
   1 & 0 & 0 & 0\\
  0 & 1 & 0 & 0\\
  0 & 0 & 1 & 0\\
  0 & 0 & 0 & -1\\
\end{pmatrix},
}
$$ it follows that it has 3-dimensional eigenspace for $+1$ and 1-dimensional eigenspace for $-1.$\par 
The fixed locus of $\alpha_{S}^{2}$ consists of $N$ curves with $a$ pairs permuted by $\alpha_{S}.$ Hence $\alpha_{S}^{*}$ on $M_{\zeta_{4}^{2}}$ has $(N-a)$-dimensional eigenspace for $+1$ and $a$-dimensional eigenspace for $-1.$
One can see that $N=k+b+2a,$ so the total effect on $h^{1,1}_{\textup{orb}}$ equals $3(N-a)+a=3k+3b+4a.$
\newline
\begin{center}
\def\arraystretch{1.4}
\begin{tabular}{c|c}
  
  Element of $C_{4}$ & $\alpha^{2}$\\
  \hline
  Irreducible components & $4N$ curves\\
  \hline
  The age & curves: 1 \\
  \hline
  Summand of $h^{1,1}_{\textup{orb}}$ & $3k+3b+4a$ \\
  
\end{tabular} 
\end{center}
\vspace{2mm}

From the orbifold formula follows that 
$$
h^{1,1}_{\textup{orb}}=r+1+4k+2(n_{1}+n_{2})+3k+3b+4a=1+r+7k+3b+2(n_{1}+n_{2})+4a.
$$

In order to compute $h^{1,2}_{\textup{orb}}$ we will use the orbifold Euler characteristic. In the table below we collect all possible intersections $X^{g}\cap X^{h},$ where $(g,h)\in C_{4}^{2}.$

\begin{center}
\def\arraystretch{1.4}
\begin{tabular}{c|c|c|c|c}
                                        & 1 & $\zeta_{4}^{}$ & $\zeta_{4}^{2}$ & $\zeta_{4}^{3}$ \\ \hline
\multicolumn{1}{c|}{1}                 & $X$ &   ${X^{\zeta_{4}^{}}}$ &  ${X^{\zeta_{4}^{2}}}$ & ${X^{\zeta_{4}^{}}}$  \\ \hline
\multicolumn{1}{c|}{$\zeta_{4}^{}$}    & ${X^{\zeta_{4}^{}}}$ &   ${X^{\zeta_{4}^{}}}$ &  ${X^{\zeta_{4}^{}}}$ & ${X^{\zeta_{4}^{}}}$  \\ \hline
\multicolumn{1}{c|}{$\zeta_{4}^{2}$}   & ${X^{\zeta_{4}^{2}}}$ & ${X^{\zeta_{4}}}$ &  ${X^{\zeta_{4}^{2}}}$ & ${X^{\zeta_{4}^{}}}$  \\ \hline
\multicolumn{1}{c|}{$\zeta_{4}^{3}$}   & ${X^{\zeta_{4}^{}}}$ &  ${X^{\zeta_{4}^{}}}$ &  ${X^{\zeta_{4}^{}}}$ & ${X^{\zeta_{4}^{}}}$ \\ 
\end{tabular}
\end{center}

By \ref{TR3} we obtain the formula
$$ e(\widetilde{X/C_{4}})=e_{\orb}(X/C_{4})=\frac{1}{4}\left(12e(X^{\zeta_{4}})+3e(X^{\zeta_{4}^{2}})\right)=6e(S^{\zeta_{4}})+3e(S^{\zeta_{4}^{2}}).$$
\par
If $D$ is of the first type, then by the Riemann-Hurwitz formula we conclude that
$$e(S^{\zeta_{4}})=2-2g(D)+2(k-1)+n_{1}\;\;\; \textup{and}\;\;\; e(S^{\zeta_{4}^{2}})=2-2g(D)+2(N-1).$$
Since $\widetilde{X/C_{4}}$ is Calabi-Yau we gets 
\begin{align*}
&h^{1,2}(\widetilde{X/C_{4}})=h^{1,1}(\widetilde{X/C_{4}})-\frac{1}{2}e(\widetilde{X/C_{4}})=1+r+7k+3b+2(n_{1}+n_{2})+4a-\\&-(-9g(D)+6k+3N+3n_{1})=1+r+k+3b+2n_{2}-n_{1}+4a+9g(D)-3N.
\end{align*}
By (\cite{AAS}, Thm. 1.1) and (\cite{AAS}, Prop. 1) we have the following relations
$$ r=\frac{1}{2}(12+k+2a+b-g(D)+4h)\;\; \textup{and}\;\; m=\frac{1}{2}(12-k-2a-b+g(D)),$$ where $\displaystyle h=\sum_{C\subseteq S^{\alpha_{S}}}(1-g(C)).$ Moreover since $D$ is of the first type $h=k-g(D),$ $n_{2}=0,$ $n_{1}=2h+4$ and $\displaystyle b=\frac{n_{1}}{2}$, thus
\begin{align*}
h^{1,2}(\widetilde{X/C_{4}})&=1+r+k+3b+2n_{2}-n_{1}+4a+9g(D)-3N=m-1+7g(D).
\end{align*}
\par
If $D$ is of the second type, then analogously the Riemann-Hurwitz formula yields
$$e(S^{\zeta_{4}})=2h+n_{1}+n_{2}\;\;\; \textup{and}\;\;\; e(S^{\zeta_{4}^{2}})=2-2g(D)+2(N-1).$$
Thus 
\begin{align*}
h^{1,2}(\widetilde{X/C_{4}})&=h^{1,1}(\widetilde{X/C_{4}})-(6h+3n_{1}+3n_{2}+3N-3g(D))=\\&=1+r+7k+3b-n_{1}-n_{2}+4a-6h-3N+3g(D).
\end{align*}
Using the additional relations $h=k,$ $n_{1}+n_{2}=2h+4$ and $b=\displaystyle \frac{n_{1}}{2}+1,$ we get 
$$h^{1,2}(\widetilde{X/C_{4}})=m+2g(D)-\frac{n_{2}}{2}.$$

Hence we proved the following theorem:
\begin{Thm}[\cite{C}, Proposition 6.3]{}
If $S^{\alpha_{S}^{2}}$ is not a sum of two elliptic curves, then for any crepant resolution of variety $X/C_{4}$ the following formulas hold
\begin{itemize}
\item If $D$ is of the first type, then 
\begin{align*}&h^{1,1}=1+r+7k+3b+2(n_{1}+n_{2})+4a,\\
&h^{1,2}=m-1+7g(D).\end{align*}
\item If $D$ is of the second type, then
\begin{align*}&h^{1,1}=1+r+7k+3b+2(n_{1}+n_{2})+4a,\\
&h^{1,2}=m+2g(D)-\frac{n_{2}}{2}.\end{align*}
\end{itemize}
\end{Thm}

\subsection{Order 6}
Let $(S, \gamma_{S})$ be a $K3$ surface with purely non-symplectic automorphism $\gamma_{S}$ of order $6.$ Consider an elliptic curve $E$ with the Weierstrass equation $y^{2}=x^{3}+1$ together with a non-symplectic automorphism $\gamma_{E}$ of order 6 such that $$\gamma_{E}(x,y)=(\zeta_{3}^{2}x, -y).$$\par

We shall keep the notation of \cite{C}.
\begin{enumerate}[label=-]
\item[] $X=S\times E,$
\item[] $P$ -- the infinity point of $E,$
\item[] $r=\dim H^{2}(S,\mathbb{C})^{\gamma_{S}},$
\item[] $m=\dim H^{2}(S,\mathbb{C})_{\zeta_{6}^{i}}$ for $i\in\{1,2,\ldots, 5\},$
\item[] $l$ -- number of curves fixed by $\gamma_{S},$
\item[] $k$ -- number of curves fixed by $\gamma_{S}^{2},$
\item[] $N$ -- number of curves fixed by $\gamma_{S}^{3},$
\item[] $p_{(2,5)}+p_{(3,4)}$ -- number of isolated points fixed by $\gamma_{S}$ of type $(2,5)$ and $(3,4)$ i.e. the action of $\gamma_{S}$ near the point linearizes to respectively
$$
{
\begin{pmatrix}
  \zeta_{6}^{2} & 0 \\
  0 & \zeta_{6}^{5} \\
\end{pmatrix}
}
\ \quad \text{and}\quad
{
\begin{pmatrix}
  \zeta_{6}^{3} & 0 \\
  0 & \zeta_{6}^{4} \\
\end{pmatrix}
},
$$ 
\item[] $n$ -- number of isolated points fixed by $\gamma_{S}^{2},$
\item[] $2n'$ -- number of isolated points fixed by $\gamma_{S}^{2}$ and switched by $\gamma_{S},$
\item[] $a$ -- number of triples $(A, A', A'')$ of curves fixed by $\gamma_{S}^{3}$ such that $\gamma_{S}(A)=A'$ and $\gamma_{S}(A')=A'',$
\item[] $b$ -- number of pairs $(B,B')$ of curves fixed by $\alpha_{S}^{2}$ such that $\gamma_{S}(B)=B',$ 
\item[] $D$ -- the curve with the highest genus in the fixed locus of $\gamma_{S},$
\item[] $G$ -- the curve with the highest genus in the fixed locus of $\gamma_{S}^{2},$
\item[] $F_{1}, F_{2}$ -- the curves with the highest genus in the fixed locus of $\gamma_{S}^{3},$
\end{enumerate}

\begin{Rem}\label{D} From (\cite{Dil}, Thm. 4.1) follows that $g(D)\in \{0,1\}.$ Moreover by \cite{MAS} we see that if $g(F_{1})\neq 0,$ $g(F_{2})\neq 0,$ then $g(F_{1})=g(F_{2})=1.$ Clearly if $g(D)=1,$ then $D= G\in \{F_{1},F_{2}\}.$ \end{Rem}

Denote by $\gamma$ an automorphism $\gamma_{S}\times \gamma_{E}^{5}$. Clearly $\langle \gamma \rangle \simeq  C_{6}.$
Similar computations as in the previous cases imply that $$h^{1,1}_{\orb}(X)^{C_{6}}=r+1\;\;\; \textrm{and} \;\;\; h^{2,1}_{\orb}(X)^{C_{6}}=m-1.$$
For any $g\in C_{6}$ let $$M_{g}:=\bigoplus_{U\in \Lambda(g)} H^{1- \age(g),\; 1-\age(g)}(U).$$ 
\par
\textsl{The action of $\gamma$ and $\gamma^{5}$.}
An automorphism $\gamma$ acts on $E$ by $\gamma(x,y)=(\zeta_{6}^{4}x, -y),$ hence it has only one fixed point --- $P$. The fixed locus of $\gamma$ on $S$ consists of $l$ curves and $p_{(2,5)}+p_{(3,4)}$ isolated points. Locally the action of $\gamma$ on $S\times E$ along the curve can be diagonalised to a matrix
$$
{
\begin{pmatrix}
   1 & 0 & 0\\
  0 & \zeta_{6} & 0 \\
  0 & 0 & \zeta_{6}^{5} \\
\end{pmatrix},
}
$$ with age equal to $1.$ In the fixed point of type $(2,5)$ and $(3,4)$ we get respectively a matrices 
$$
{
\begin{pmatrix}
   \zeta_{6}^{2} & 0 & 0\\
  0 & \zeta_{6}^{5} & 0 \\
  0 & 0 & \zeta_{6}^{5} \\
\end{pmatrix}
}
\;
\textup{ and }
\;
{
\begin{pmatrix}
   \zeta_{6}^{3} & 0 & 0\\
  0 & \zeta_{6}^{4} & 0 \\
  0 & 0 & \zeta_{6}^{5} \\
\end{pmatrix},
}
$$
hence their ages equal 2. 
\par In case of the action of $\gamma^{5}$ we observe that locus consists of $l$ curves and $p_{(2,5)}+p_{(3,4)}$ isolated points. Along the curve we have a matrix 
$$
{
\begin{pmatrix}
   1 & 0 & 0\\
  0 & \zeta_{6}^{5} & 0 \\
  0 & 0 & \zeta_{6} \\
\end{pmatrix},
}
$$ while in fixed points we have a matrices
$$
{
\begin{pmatrix}
   \zeta_{6}^{1} & 0 & 0\\
  0 & \zeta_{6}^{4}  & 0 \\
  0 & 0 & \zeta_{6} \\
\end{pmatrix}
}
\ \quad \textrm{or}\quad
{
\begin{pmatrix}
   \zeta_{6}^{2} & 0 & 0\\
  0 & \zeta_{6}^{3}  & 0 \\
  0 & 0 & \zeta_{6} \\
\end{pmatrix},
}
$$ 
hence their ages are equal to 1.
The above analysis shows that the effect on $h^{1,1}_{\textup{orb}}$ from both actions equals $l+l+p_{(2,5)}+p_{(3,4)}=2l+p_{(2,5)}+p_{(3,4)}$.
\newline
\begin{center}
\def\arraystretch{1.4}
\begin{tabular}{c|c|c}
  
  Element of $C_{6}$ & $\gamma$ & $\gamma^{5}$\\
  \hline
  Irr. comp. & $3l$ curves, $3p_{(2,5)}+3p_{(3,4)}$ pts. & $3l$ curves, $3p_{(2,5)}+3p_{(3,4)}$ pts.\\
  \hline
  The age & curve: 1, point: 2 & curve: 1, point: 1\\
  \hline
  Summand of $h^{1,1}_{\textup{orb}}$ & $l$ & $l+p_{(2,5)}+p_{(3,4)}$\\
  
\end{tabular} 
\end{center}
\par
\textsl{The action of $\gamma^{2}$ and $\gamma^{4}$.} Automorphisms $\gamma^{2}$ and $\gamma^{4}$ act on $E$ by $\gamma^{2}(x,y)=(\zeta_{6}^{2}x, y)$ and  $\gamma^{4}(x,y)=(\zeta_{6}^{4}x, y)$, they have three fixed points --- $\{P, (0,i), (0,-i)\}$ from which only $P$ is invariant under $\gamma$ and the remaining two are switched. Identifying $M_{\zeta_{6}^{2}}$ and $M_{\zeta_{6}^{4}}$ with the vector space spanning by irreducible components we will find the action of induced map $\gamma^{*}$ on it. \par The matrix of the action of $\gamma^{*}$ on $M_{\zeta_{6}^{2}}$ is
$$
{
\begin{pmatrix}
   1 & 0 & 0\\
  0 & 0 & 1 \\
  0 & 1 & 0 \\
\end{pmatrix}\sim 
}
{
\begin{pmatrix}
   1 & 0 & 0\\
  0 & 1 & 0 \\
  0 & 0 & -1 \\
\end{pmatrix}
}
$$ hence it produces a 2-dimensional eigenspace for $+1$ and 1-dimensional eigenspace for $-1.$ We have similar decomposition in the case of the action $\gamma^{*}$ on $M_{\zeta_{6}^{2}}$. \par The fixed locus of $\gamma_{S}^{2}$ consists of $p_{(2,5)}$ and $n'$ pairs of points switched by $\gamma_{S}.$ Locally, the action of $\gamma^{2}$ at point has a matrix 
$${
\begin{pmatrix}
   \zeta_{6}^{4} & 0 & 0\\
  0 & \zeta_{6}^{4}  & 0 \\
  0 & 0 & \zeta_{6}^{4} \\
\end{pmatrix},
}
$$ with age equals 2. In the case of automorphism $\gamma^{4}$ we have a matrix
$$
{
\begin{pmatrix}
   \zeta_{6}^{2} & 0 & 0\\
  0 & \zeta_{6}^{2}  & 0 \\
  0 & 0 & \zeta_{6}^{2} \\
\end{pmatrix},
}
$$ with age equal 1. \par
Notice that both $\gamma_{S}^{2}$ and $\gamma_{S}^{4}$ has the same fixed points but with different ages. Thus $\gamma_{S}^{*}$ on $M_{\zeta_{6}^{2}}$ has $(p_{2,5}^{1}+n_{1})$-dimensional eigenspace for $+1$ and $n'_{1}$-dimensional eigenspace for $-1,$ while $\alpha_{S}^{*}$ on $M_{\zeta_{6}^{4}}$ has $(p_{2,5}^{2}+n_{2})$-dimensional eigenspace for $+1$ and $n_{2}$-dimensional eigenspace for $-1,$ where $p_{2,5}^{1},$ $p_{2,5}^{2},$ $n'_{1}$ and $n'_{2}$ are naturals defined as $n'_{1}+n'_{2}=n'$ and $p_{2,5}^{1}+p_{2,5}^{2}=p_{2,5}.$
\par 
Moreover, $a$ pairs of curves in the locus of $\gamma_{S}^{2}$ are switched by $\gamma.$ Hence $\gamma_{S}^{*}$ has $(k-b)$-dimensional eigenspace for $+1$ and $b$-dimensional eigenspace for $-1.$ The same decomposition we will obtain in case of the action of $\gamma^{4}.$
\par
By the K{\"u}nneth formula both actions add to $h^{1,1}_{\orb}$ \begin{align*}&2(2(k-b)+b)+ 2(p_{2,5}^{1}+n'_{1})+n_{1}'+2(p_{2,5}^{2}+n'_{2})+n_{2}' = \\&=4k-2b+2p_{(2,5)}+3n'.\end{align*}

\begin{center}
\def\arraystretch{1.4}
\begin{tabular}{c|c|ccc}
Element of $C_{6}$     & $\gamma^{2}$ & $\gamma^{4}$ &  \\ \cline{1-3}
Irreducible components & $3k$ curves, $3n$ points & $3k$ curves, $3n$ points &  \\ \cline{1-3}
The age & curve: 1, point: 2 & curve: 1, point: 1 &  \\ \cline{1-3}
Summand of $h^{1,1}_{\orb}$ &  \multicolumn{2}{c}{$4k-2b+2p_{(2,5)}+3n'$}  &  &  \\ 
\end{tabular}
\end{center}
\vspace{2mm}
\par
\textsl{The action of $\gamma^{3}.$} An automorphism $\gamma^{3}$ acts on $E$ by $\gamma^{3}(x,y)=(x, -y),$ hence it has four fixed points --- $\{P, (1,0), (\zeta_{3},0), (\zeta_{3}^{2},0)\}$ from which only $P$ is invariant under $\gamma$ and the remaining three form a 3-cycle. Thus the action of ${\gamma_{E}^{5}}^{*}$ on $M_{\zeta_{6}^{3}}$ has the matrix 
$$
{
\begin{pmatrix}
   1 & 0 & 0 & 0\\
  0 & 0& 0 & 1 \\
  0 &1 & 0 &0 \\
  0 &0 & 1 &0 \\
\end{pmatrix}\sim 
}
{
\begin{pmatrix}
   1 & 0 & 0 & 0\\
  0 & 1& 0 & 0 \\
  0 &0 & \zeta_{3} &0 \\
  0 &0 & 0 &\zeta_{3}^{2} \\
\end{pmatrix}
}
$$ so it produces 2-dimensional eigenspace for 1, 1-dimensional eigenspace for $\zeta_{3}$ and 1-dimensional eigenspace for $\zeta_{3}^{2}.$ \par
Clearly the local action of $\gamma^{3}$ on curve has age equal to 1. The locus of $\gamma_{S}^{3}$ consists of $N$ curves with $a$ triples of curves permuted by $\gamma_{S}.$ Thus $\gamma_{S}^{*}$ on $M_{\zeta_{6}^{3}}$ has $(N-2a)$-dimensional eigenspace for +1 and two $a$-dimensional eigenspaces for $\zeta_{3}$ and $\zeta_{3}^{2}.$ 
\par
The effect on $h^{1,1}_{\orb}$ equals $2(N-2a)+a+a=2N-2a.$
\vspace{2mm}
\begin{center}
\def\arraystretch{1.4}
\begin{tabular}{c|c}
  Element of $C_{6}$ & $\gamma^{3}$\\
  \hline
  Irreducible components & $4N$ curves\\
  \hline
  The age & curve: 1 \\
  \hline
  Summand of $h^{1,1}_{\orb}$ & $2N-2a$ \\
\end{tabular} 
\end{center}

Consequently, by the orbifold formula we see that
\begin{align*}
h^{1,1}_{\orb}&=r+1+l+l+p_{(2,5)}+p_{(3,4)}+4k-2b+2p_{(2,5)}+3n'+2N-2a=\\&=r+1+2l+2N-2b+4k-2a+3n'+3p_{(2,5)}+p_{(3,4)}.
\end{align*}

By the orbifold cohomology formula we see that non-zero contribution to $h^{1,2}_{\textup{orb}}$ have only curves in $\Lambda(g)$ for any $g\in C_{6}.$ \par If $g(D)\geq 1,$ by \ref{D} we can assume that $D=G=F_{1}.$ We see that contributions of $\gamma$ and $\gamma^{5}$ equal to $2g(D).$ The automorphisms $\gamma^{2}$ and $\gamma^{4}$ have three fixed points on $E$ with one $2$-point orbit, hence by K{\"u}nneth's formula the effect on $h^{1,2}_{\textup{orb}}$ in this case equals $4g(D).$ Since $\gamma^{3}$ has four fixed points with 3-points orbit, we find that it's contribution is equal to $2g(D)+g(F_{2})+g(F_{2}/\gamma_{S}).$ Thus $$h^{1,2}_{\textup{orb}}(\widetilde{X/C_{6}})=m-1+8g(D)+g(F_{2})+g(F_{2}/\gamma_{S}).$$ \par
Now consider the case $g(D)=0.$ By the same argument as above, we see that contribution of $\gamma^{2}$ and $\gamma^{4}$ equals $2g(G)+2g(G/\gamma_{S}),$ while the summand from $\gamma^{3}$ equals $g(F_{21})+g(F_{1}/\gamma_{S})+g(F_{2})+g(F_{2}/\gamma_{S}),$ hence
$$h^{1,2}_{\textup{orb}}(\widetilde{X/C_{6}})=m-1+2g(G)+2g(G/\gamma_{S})+g(F_{1})+g(F_{1}/\gamma_{S})+g(F_{2})+g(F_{2}/\gamma_{S}).$$\par

We proved the following theorem:

\begin{Thm}[\cite{C}, Proposition 7.3]{}
For any crepant resolution of variety $X/C_{6}$ the following formulas hold 
\begin{align*}&h^{1,1}=r+1+2l+2N-2b+4k-2a+3n'+3p_{(2,5)}+p_{(3,4)},\\
&h^{1,2}=\begin{cases} m-1+8g(D)+g(F_{2})+g(F_{2}/\gamma_{S})\;\; &\textrm{if } g(D)\geq 1, \\ m-1+2g(G)+2g(G/\gamma_{S})+g(F_{1})+g(F_{1}/\gamma_{S})\\+g(F_{2})+g(F_{2}/\gamma_{S})\;\; &\textrm{if } g(D)=0.\end{cases} \end{align*}
\end{Thm}

Now, we will compute $h^{2,1}_{\orb}$ using orbifold Euler characteristic.
In the table below we collect all possible intersections $X^{g}\cap X^{h},$ where $(g,h)\in C_{6}^{2}.$

\begin{center}
\def\arraystretch{1.4}
\begin{tabular}{c|c|c|c|c|c|c}
        & 1 & $\zeta_{6}$ & $\zeta_{6}^{2}$ &$\zeta_{6}^{3}$& $\zeta_{6}^{4}$ & $\zeta_{6}^{5}$ \\ \hline
\multicolumn{1}{c|}{1}                                  & $X$  &{$X^{\zeta_{6}}$}&           {$X^{\zeta_{6}^{2}}$ }                      &  {$X^{\zeta_{6}^{3}} $}                          &                  {$X^{\zeta_{6}^{2}}$ }                 &           {$X^{\zeta_{6}}$}               \\ \hline
\multicolumn{1}{c|}{$\zeta_{6}$}                       &{$X^{\zeta_{6}}$}& {$X^{\zeta_{6}}$}        &    {$X^{\zeta_{6}}$} &          {$X^{\zeta_{6}}$}                  &       {$X^{\zeta_{6}}$}                       &        {$X^{\zeta_{6}}$}   \\ \hline
\multicolumn{1}{c|}{$\zeta_{6}^{2}$} &  {$X^{\zeta_{6}^{2}}$ }   &   {$X^{\zeta_{6}}$}      &     {$X^{\zeta_{6}^{2}}$ }           &         {$X^{\zeta_{6}}$}                   &      {$X^{\zeta_{6}^{2}}$ }                    & {$X^{\zeta_{6}}$} \\ \hline
\multicolumn{1}{c|}{$\zeta_{6}^{3}$} & {$X^{\zeta_{6}^{3}} $} &  {$X^{\zeta_{6}}$}        &{$X^{\zeta_{6}}$}&    {$X^{\zeta_{6}^{3}} $}                           &      {$X^{\zeta_{6}}$}                     &       {$X^{\zeta_{6}}$}              \\ \hline
\multicolumn{1}{c|}{$\zeta_{6}^{4}$} & {$X^{\zeta_{6}^{2}}$ }   &  {$X^{\zeta_{6}}$}        & {$X^{\zeta_{6}^{2}}$ } &   {$X^{\zeta_{6}}$}                         &   {$X^{\zeta_{6}}$}                            &       {$X^{\zeta_{6}}$} \\ \hline
\multicolumn{1}{c|}{$\zeta_{6}^{5}$} & {$X^{\zeta_{6}}$}  &  {$X^{\zeta_{6}}$}           &{$X^{\zeta_{6}}$}&{$X^{\zeta_{6}}$} & {$X^{\zeta_{6}}$}  &{$X^{\zeta_{6}}$}    \\ 
\end{tabular}
\end{center}

From \ref{TR3} we see that 
$$ e(\widetilde{X/C_{6}})=\frac{1}{6}\left(24e(X^{\zeta_{6}^{}})+8e(X^{\zeta_{6}^{2}})+3e(X^{\zeta_{6}^{3}})\right)=4e(S^{\zeta_{6}})+4e(S^{\zeta_{6}^{2}})+2e(S^{\zeta_{6}^{3}}).$$
Thus $$h^{1,2}(\widetilde{X/C_{6}})=h^{1,1}(\widetilde{X/C_{6}})-2e(S^{\zeta_{6}})-2e(S^{\zeta_{6}^{2}})-e(S^{\zeta_{6}^{3}}).$$
\par
By the Riemann-Hurwitz formula we obtain \begin{align*} &e(S^{\zeta_{6}})=2(l-1)+2-2g(D)+p_{(2,5)}+p_{(3,4)},\\ &e(S^{\zeta_{6}^{2}})=2(k-1)+1-g(G)+1-g(G/\gamma_{S})+n,\\ &e(S^{\zeta_{6}^{3}})=2(N-2)+1-g(F_{1})+1-g(F_{1}/\gamma_{s})+1-g(F_{2})+1-g(F_{2}/\gamma_{s})\end{align*}
hence after simplifying 
\begin{align*} h^{1,2}(\widetilde{X/C_{6}})&=r+1-2l-2b-2a+3n'+p_{(2,5)}-p_{(3,4)}-2n+4g(D)+2g(G)+\\&+2g(G/\gamma_{s})+g(F_{1})+g(F_{1}/\gamma_{s})+g(F_{2})+g(F_{2}/\gamma_{s}).\end{align*}

Comparing both formulas we get:

\begin{Cor} With the notation above, the following relation holds
 $$-m+r+2-2l-2b-2a+3n'+p_{(2,5)}-p_{(3,4)}-2n+4g(D)=0.$$ 
\end{Cor}

\subsection*{Acknowledgments} This paper is a part of author's master thesis. I am deeply grateful to my advisor S\l{}awomir Cynk for recommending me to learn this area and his help.

\end{document}